\documentclass[11pt]{article}

\usepackage{latexsym}
\usepackage{amsfonts}
\usepackage{amsmath}
\usepackage{amsthm}
\usepackage{amssymb}
\usepackage{mathrsfs}
\usepackage{dsfont}
\usepackage[english]{babel}
\usepackage{caption}
\usepackage{epsfig}
\usepackage{float}
\usepackage{graphicx}
\usepackage{epsfig}
\usepackage{amsbsy}
\usepackage{hyperref}
\usepackage{subfig}
\usepackage{graphicx,color}

\newtheorem*{prop*}{Theorem}
\newtheorem{theo}{Theorem}

\newtheorem{hyp}{Hypothesis}

\newcommand{\zerarcounters}{\setcounter{equation}{0}\setcounter{theorem}{0}}

\newcommand{\beq}{\begin{equation}}
\newcommand{\eeq}{\end{equation}}

\newcommand{\ZZZ}{\mathds{Z}}

\newcommand{\RRR}{\mathds{R}}
\newcommand{\TTT}{\mathds{T}}

\newcommand{\calA}{{\mathcal A}}
\newcommand{\BB}{{\mathcal B}}
\newcommand{\CCCC}{{\mathcal C}}
\newcommand{\calD}{{\mathcal D}}

\newcommand{\calH}{{\mathcal H}}

\newcommand{\gotA}{{\mathfrak A}}

\newcommand{\gotD}{{\mathfrak D}}

\newcommand{\io}{\infty}
\newcommand{\e}{\varepsilon}
\newcommand{\al}{\alpha}

\newcommand{\m}{\mu}

\newcommand{\p}{\pi}
\newcommand{\ka}{\kappa}

\newcommand{\om}{\omega}

\newcommand{\f}{\varphi}

\def\ins#1#2#3{\vbox to0pt{\kern-#2 \hbox{\kern#1 #3}\vss}\nointerlineskip}

\begin{document}

\title{\bf Stable dynamics in forced systems with sufficiently high/low forcing frequency}

\author{{\bf M. Bartuccelli$^{1}$, G. Gentile$^{2}$, J.A. Wright$^{1}$}
\\ \small
$^{1}$ Department of Mathematics, University of Surrey, Guildford, GU2 7XH, UK
\\ \small
$^{2}$ Dipartimento di Matematica e Fisica, Universit\`a Roma Tre, Roma, 00148, Italy
\\ \small
E-mail:  m.bartuccelli@surrey.ac.uk, gentile@mat.uniroma3.it, j.wright@surrey.ac.uk}

\date{}

\maketitle

\begin{abstract}
We consider a class of parametrically forced Hamiltonian systems with one-and-a-half degrees
of freedom and study the stability of the dynamics when the frequency of the forcing is
relatively high or low. We show that, provided the frequency of the forcing is sufficiently high,
KAM theorem may be applied even when the forcing amplitude is far away from the perturbation regime.
A similar result is obtained for sufficiently low frequency forcing, but in that case
we need the amplitude of the forcing to be not too large; however we are still able to
consider amplitudes of the forcing which are outside of the perturbation regime.
Our results are illustrated by means of numerical simulations for the system of a forced cubic oscillator.
In addition, we find numerically that the dynamics are stable even when the forcing amplitude is very large
(beyond the range of validity of the analytical results), provided the frequency of the forcing
is taken correspondingly low.
\end{abstract}

\zerarcounters
\section{Introduction}
\label{sec:1}

Forced systems can display unexpected behaviour in the extreme cases in which the oscillation
of the forcing is either very slow or very fast. It is well known that, in such situations,
the dynamics can be strikingly different with respect to the case in which the forcing period is
comparable with that of the unperturbed motions \cite{BM1,PR,C,AKN,B2,Ch}.
A classical example is provided by the pendulum with oscillating support \cite{Ste,B1,K,BM1,BM2,BGG}.
The standard technique used to attack the problem is the averaging method; however, a rigorous
implementation is quite not trivial, as it involves dealing with small divisor problems.
Indeed, the main idea underlying the averaging is that the forced system, in suitable coordinates,
can be considered as a perturbation of an integrable system, so that KAM-like arguments apply \cite{N2,AKN}.

For a periodically forced system, if the amplitude of the forcing is small, then the KAM theorem yields that most
of the unperturbed tori persist. On the basis of heuristic arguments, one expects that something of the same kind still occurs
even when the amplitude of the forcing is not small, provided its frequency is large or small enough.
Such an expectation is ultimately based on averaging arguments. However, to reduce the problem to results available
in the literature requires some analysis. Since this has not been done explicitly, as far as we are aware of,
we think it may be of interest to discuss how the argument works.

The case of large frequency can be reduced quite easily to Neishdatd's averaging theorem. However,
at a certain point, we shall need an additional close-to-identity transformation which strongly relies
one the one-dimensionality of the system (more precisely on its integrability).
The case of low frequency is related to the problem of boundedness of the solutions
of forced systems in phase space \cite{DZ2,Le1,LL,Le2,LZ}. However, to prove boundedness
one has to prove the existence of confining KAM tori far away from the origin. Hence only the asymptotic behaviour
of the potential really counts, and in general one needs a condition on the growth of the potential at infinity -- besides
smoothness conditions, see \cite{Li,Le2}. On the contrary, to study the existence of KAM tori in a fixed
region of phase space of a forced system, one needs information about the potential in that region: this explains
why the assumptions we shall require on the potential are stronger, as they are not just asymptotic properties.

For clarity purposes, in this paper we focus on a case study, the forced cubic oscillator,
which has been extensively investigated in the literature all along as a paradigmatic model \cite{M,BBDGG,BDG,C}.
We shall prove that most of the phase space is filled by KAM invariant tori when the period of the forcing is sufficiently
large or sufficiently small. In particular we find that this happens for values of the forcing amplitude far beyond the perturbation regime.
Furthermore, we provide numerical results which strongly suggest that the amplitude can be allowed to be still larger
than the analytical bounds. We leave as an open problem how to improve the analytical bounds,
so as to give full explanation of the numerical findings.
We briefly discuss at the end how to relax the assumptions on the potential and extend the results to more general systems.

A cubic oscillator subject to a periodic driving force is described by the Hamiltonian
\begin{equation} \label{eq:1.1}
H(y,x,t)= \frac{y^2}{2} + \left( 1 + \m f(\omega t) \right) \frac{x^{4} }{4} ,
\end{equation}
where $(y,x)\in\RRR^{2}$, $t$ denotes the time, the \emph{driving force} $f$ is a $2\p$-periodic analytic function
of its argument, with zero average and $\|f\|_{\io}=1$; $\mu\in \RRR$ and $\omega\in\RRR$ are, respectively,
the \emph{amplitude} and the \emph{frequency} of the driving force.
The phase space for the system is $\RRR^{2} \times \TTT_{\om}$,
with $\TTT_{\om}=\RRR/ (2\pi/\om)\ZZZ$ and the corresponding equations of motion are
\begin{equation} \label{eq:1.2}
\begin{cases}
\dot x = y , \\
\dot y = - \left( 1 + \m f(\omega t) \right) x^{3} , \\
\dot t = 1 ,
\end{cases}
\end{equation}
where the dot denotes derivative with respect to $t$.

For $\mu=0$ the system is integrable, so that the full phase space is filled by invariant tori,
with $t$ increasing linearly in time with frequency $1$ and $(x,y)$ moving on a closed orbit $\CCCC$
in the plane with frequency $\Omega$ depending on the initial data.

Fix arbitrarily $\Omega_{2}>\Omega_{1}>0$ and consider, for the unperturbed cubic oscillator, the closed,
concentric orbits
$\CCCC_{1}$ and $\CCCC_{2}$ which run with frequencies $\Omega_1$ and $\Omega_{2}$, respectively.
Call $\calD_{0} \subset \RRR^{2}$ the bounded region enclosed between $\CCCC_{1}$ and $\CCCC_{2}$,
and set $\calD=\calD_{0} \times \TTT_{\om}$. The reasons behind the choice of the domain $\calD$
are the following: we want to fix a bounded region in phase space so as to estimate the relative measure
of the persisting tori (the curve $\CCCC_{2}$ is a natural boundary for such a region) and at any rate
one is forced to exclude a small region around the origin, where chaotic dynamics may become dominant
(the role of the internal curve $\CCCC_{1}$ is just to cut off such a small region since the beginning).

For $\mu$ small enough, we can apply KAM theorem \cite{AKN},
so as to obtain the persistence of most of the invariant tori inside $\calD$,
independently of the value of the frequency $\om$ (and of the value of $\Omega_1$ as well:
$\Omega_1$ can be taken arbitrarily small).

\begin{theo} \label{thm:1}
For $\mu$ small enough the set of persisting invariant tori in $\calD$ for the system with Hamiltonian \eqref{eq:1.1}
leaves out a set with relative measure $O(\sqrt{\mu})$.
\end{theo}

In this note we want to show that the similar results still hold when removing the condition
of smallness on $\mu$, provided the frequency $\om$ is either large enough or small enough.
This will be discussed, respectively, in Section \ref{sec:2} and Section \ref{sec:3}.
In particular, by increasing $\omega$ from $0$ to infinity, a double transition regularity-chaos-regularity is
expected to occur. We support such an expectation by providing in Section \ref{sec:4} numerical results,
which also give evidence that the regularity regime extends to wider ranges
of the parameters for which the analytical results do not apply. Finally in Section \ref{sec:5}
we discuss how to extend our analysis to more general systems.

\zerarcounters
\section{High frequency regime}
\label{sec:2}

We consider first the case of $\om$ large and we set $\om=1/\e$, with $\e$ small.
We can formulate our result for such a case in a more general context.
In action-angle variables, the Hamiltonian \eqref{eq:1.1} becomes
\begin{equation} \label{eq:2.1}
H(A,\al,t)= H_{0}(A) + \m H_{1}(A,\al,\omega t) ,
\end{equation}
where
\begin{equation} \label{eq:2.2}
H_{0}(A)= \frac{1}{4} \left( \frac{3A}{T} \right)^{4/3} , \qquad
H_{1}(A,\al,t) = \frac{1}{4} \left( \frac{3A}{T} \right)^{4/3} {\rm cn}^{4}(T\al) \, f(t) ,
\end{equation}
with ${\rm cn}\,\al :={\rm cn}(\al,1/\sqrt{2})$ and $T:=4K(1/\sqrt{2})/2\pi$,  where ${\rm cn}(\al,k)$ and
$K(k)$ denote respectively the cosine-amplitude function and the complete elliptic integral of the first kind
with elliptic modulus $k$.

The corresponding Hamilton equations are
\begin{equation} \label{eq:2.3}
\begin{cases}
\dot \al = \Omega_0(A) + \mu \partial_A H_{1}(A,\al,\om t) , \\
\dot A = - \mu \partial_\al H_{1}(A,\al,\om t) ,
\end{cases}
\end{equation}
where $\Omega_0(A):=\partial_{A} H_0(A)=(3A/T^4)^{1/3}$.
In terms of the variables $(A,\al,t)$ the domain $\calD$ defined in Section \ref{sec:1} becomes
$\gotD=\gotA_{0} \times \TTT \times \TTT_{\om}$, where $\TTT=\TTT_{1}$ and $\gotA_{0}=\{A\in \RRR_{+} :
A_{1} \le |A| \le A_{2}\}$, with $A_{1}=\Omega_{1}^{3}T^4/3$ and $A_{2}=\Omega_{2}^{3}T^4/3$.

More generally we shall consider Hamiltonians of the form \eqref{eq:2.1} and make
the following assumptions on $H_{0}$ and $H_{1}$ --- trivially satisfied in the case \eqref{eq:2.2}.

\begin{hyp} \label{hyp:1}
Assume the Hamiltonian function \eqref{eq:2.1} to be real-analytic in a domain
$\gotD:=\gotA_{0} \times \TTT \times \TTT_{\om}$, where $\gotA_{0}$ is an open subset of $\RRR$.
\end{hyp}

\begin{hyp} \label{hyp:2}
Assume $A \mapsto \Omega_0(A) := \partial_{A} H_{0}(A)$ to be a local diffeomorphism on $\gotA_{0}$.
\end{hyp}

\begin{hyp} \label{hyp:3}
Assume $\displaystyle{\langle H_{1}(A,\al,\cdot) \rangle := \int_{\TTT} \frac{{\rm d}t}{2\p} \, H_{1}(A,\al,t)  =0}$,
$\forall (A,\al)\in \gotA_{0} \times \TTT$.
\end{hyp}

In \eqref{eq:2.1} rescale time $t \to \tau=\omega t$. Then the equations of motion \eqref{eq:2.3} become
\begin{equation} \label{eq:2.4}
\begin{cases}
\dot \al = \e \Omega_0(A) + \e \mu \partial_A H_{1}(A,\al,\tau) , \\
\dot A = - \e \mu \partial_\al H_{1}(A,\al,\tau) ,
\end{cases}
\end{equation}
where now the dot denotes derivative with respect to time $\tau$.

We can apply Neishtadt's averaging theorem \cite{N1,N2} to cast the system into the form
\begin{equation} \label{eq:2.5}
\begin{cases}
\dot \al = \e \left( \Omega_0(A) + \mu \partial_A V_{\e} (A,\al) + \mu \partial_{A} R_{\e} (A,\al,\tau) \right) , \\
\dot A = - \e \left( \mu \partial_\al V_{\e} (A,\al) + \mu \partial_{\al} R_{\e} (A,\al,\tau) \right) , \\
\end{cases}
\end{equation}
where $V_{\e}$ and $R_{\e}$ are suitable analytic functions, with
$V_{\e}(A,\al)=\langle H_{1}(A,\al,\cdot) \rangle + O(\e)$ and
$R_{\e}$ an exponentially small remainder, that is
$|R_{\e}| \le C \exp (-c/\e)$ for some positive constants $c$ and $C$.
The change of coordinates is canonical and $\e$-close to the identity.
In order not to overwhelm the notations,
we denote the new variables with the same letters as the old ones.

By Hypothesis \ref{hyp:3} the average of $H_1$ vanishes, hence $V_{\e}$ is a correction
of order $\e$ to $H_{0}$. So one can perform a further close-to-identity change of coordinates
which leads to the equations
\begin{equation} \label{eq:2.6}
\begin{cases}
\dot \al = \e \overline{\Omega}_\e(A) + \e \mu \partial_{A} \overline{R}_{\e} (A,\al,\tau) , \\
\dot A = - \e \mu \partial_{\al} \overline{R}_{\e} (A,\al,\tau) , \\
\end{cases}
\end{equation}
where $\overline{\Omega}_{\e}(A)=\Omega_{0}(A)+O(\e)$ and $\overline{R}$ is still exponentially small
(and again we still denote by $(A,\al)$ the transformed coordinates). The corresponding Hamiltonian is
\begin{equation} \label{eq:2.7}
\overline{H}(A,\al,\e) = \e \left( \overline{H}_\e(A) + \mu \overline{R}_\e(A,\al,\tau) \right) .
\end{equation}
The overall change of coordinates leading to \eqref{eq:2.7} is close to the identity within $O(\e)$
and hence, up to a region with measure $O(\varepsilon)$, the domain $\gotD$ is transformed into
a region enclosed between two KAM invariant tori. By studying the Hamiltonian $\overline{H}(A,\al,\e)/\e$
we see that we can apply once more KAM theorem and conclude that most of the unperturbed tori
for the Hamiltonian $\overline{H}_\e(A)$ persist when the perturbation $\overline{R}_\e(A,\al,\tau)$
is switched on. Since now the perturbation is exponentially small, the relative measure of the tori
which are destroyed is exponentially small in $\e$. To go back to the original coordinates,\
we have to scale back time. So we obtain the following result.

\begin{theo} \label{thm:2}
Consider the system with Hamiltonian \eqref{eq:2.1} and assume Hypotheses \ref{hyp:1} to \ref{hyp:3}.
For any value of $\mu$, for $\om$ large enough the domain $\gotD$ is filled by KAM invariant tori,
up to a region whose relative measure is $O(1/\om)$. Apart from a thin region close to the boundary,
the invariant tori leave out a region with measure exponentially small in $1/\om$.
\end{theo}

Note that Hypothesis \ref{hyp:3} is needed here, contrary to Theorem \ref{thm:1}, to ensure that the averaged system
is integrable (such a condition is automatically satisfied for $\m$ small, without any assumption on $H_1$).

\zerarcounters
\section{Low frequency regime}
\label{sec:3}

Now consider \eqref{eq:1.1} with $\om=\e$. We can reason as in \cite{LL} (see also \cite{DZ1,DZ2,Le2,CNY}).
Fix $|\mu|<1$, so that $1+ \m f(\om t) >0$. We rewrite the equations of motion \eqref{eq:1.2} as
\begin{equation} \label{eq:3.1}
\begin{cases}
\dot x = y , \\
\dot y = - a(\e t) \, x^{3} , \qquad
a(t) := 1 + \m f(t) .
\end{cases}
\end{equation}
Then the argument proceeds through the following steps.

First, through a time-dependent canonical change of coordinates $(x,y) \mapsto (A,\al)$, with
\begin{equation} \label{eq:3.2}
x = \left( \frac{3A}{T} \right)^{1/3} \!\!\!\! (a(\e t))^{-1/6} {\rm cn}\,(T\al) , \qquad
y = - \left( \frac{3A}{T} \right)^{2/3} \!\!\!\! (a(\e t))^{1/6} {\rm sn}\,(T\al) \, {\rm dn} \, (T\al) ,
\end{equation}
where $\rm sn(\al)$ and ${\rm dn}(\al)$ are the sine-amplitude and delta-amplitude functions
with modulus $k=1/\sqrt{2}$, respectively,
one writes \eqref{eq:3.1} as the Hamilton equations corresponding to the Hamiltonian
\begin{subequations} \label{eq:3.3}
\begin{align}
\calH(A,\al,t) & = \calH_{0}(A,\al,t) + \frac{\e}{6} \frac{3A}{T}
\, {\rm cn} \, (T\al) \, {\rm sn} \, (T\al) \, {\rm dn} \, (T\al) \, \frac{b (\e t)}{a(\e t)} , \\
\calH_{0}(A,\al,t) & = \frac{1}{4} \left( \frac{3A}{T} \right)^{4/3} \!\!\!\! \left( a(\e t) \right)^{1/3} .
\end{align}
\end{subequations}
where $b(t)=\dot a(t)$ and $T:=4K(1/\sqrt{2})/2\pi$.
Next, in order to eliminate the dependence on time in the leading term one makes the change
of coordinates $(A,\al,t) \mapsto (p,q,s)$, with
\begin{equation} \label{eq:3.4}
p = \calH(A,\al,t) , \qquad q = t , \qquad s = \alpha ,
\end{equation}
which leads to the Hamiltonian
\begin{subequations} \label{eq:3.5}
\begin{align}
\calA(p, q , s) & = \calA_{0}(p,\e q) + \e \calA_{1}(p,\e q,s) , \\
\calA_{0}(p,q) & = \frac{T}{3} \frac{(4p)^{3/4}}{(a(q))^{1/4}} ,
\end{align}
\end{subequations}
for some function $\calA_{1}$ of order 1 in $\e$.
Then one may perform a further change of variables $(p,q,s) \mapsto (J,\f,s)$
into action-angle variables for $\calA_0$, so yielding the Hamiltonian
\begin{subequations} \label{eq:3.6}
\begin{align}
\BB(J,\f,s) & = \BB_{0}(J) + \e \BB_{1}(J,\e\f,s) , \\
\BB_{0}(J) & = \left( \ka J \right)^{3/4} , \qquad
\frac{1}{\ka} := \frac{1}{4} \left( \frac{3}{T} \right)^{4/3}
\int_{0}^{2\pi} (a(q))^{1/3} \frac{{\rm d}q}{2\pi} ,
\end{align}
\end{subequations}
for a suitable function $\BB_{1}$ of order 1 in $\e$.

Finally we integrate the Hamiltonian equations corresponding to \eqref{eq:3.6}
between $s=0$ and $s=2\pi$. Denote by $(J(s),\f(s))$ the solution;
then, defining $\psi(s)=\e \f(s)$ and setting $(J',\psi')=(J(2\pi),\psi(2\pi))$ and $(J,\psi)=(J(0),\psi(0))$,
we obtain the twist map
\begin{equation} \label{eq:3.7}
\begin{cases}
J' = J + \e^{2} F(J,\psi) , \\
\psi'=\psi + \e \Omega_{0}(J) + \e^{2} G(J,\psi) , \qquad
\displaystyle{\Omega_0(J)=\frac{3\pi}{2} \ka^{3/4} J^{-1/4}} ,
\end{cases}
\end{equation}
for suitable analytic functions $F$ and $G$.

Therefore we can apply Moser's twist theorem and conclude that any invariant
curve with Diophantine rotation number persists for $\e$ small enough.
For fixed $\e$, the relative measure of the persisting curves
in a given region of the cylinder is $O(\sqrt{\e})$ \cite{La,Sva}.
Note that, in order to prove just boundedness of the solutions,
it would be enough to prove the existence of an invariant circle of constant type
(as sometimes done in the literature, see for instance \cite{CNY}; see also the comments in \cite{LL}).
On the contrary, to prove that the persisting tori have large measure for $\e$ small,
a milder Diophantine condition is required; one could even allow a Bryuno condition on
the rotation number \cite{B3}, as done in \cite{G}, but this would not increase appreciably
the measure of the invariant curves.

Coming back to the original coordinates we obtain the existence of a large measure set
of invariant tori for the continuous flow which has the twist map \eqref{eq:3.7} as Poincar\'e section
at times multiples of $2\p$. We can summarise the discussion by the following statement.

\begin{theo} \label{thm:3}
Consider the system with Hamiltonian \eqref{eq:1.1} and fix $\mu\in(-1,1)$.
For $\om$ small enough the domain $\gotD$ is filled by KAM invariant tori,
up to a region whose relative measure is $O(\sqrt{\om})$.
\end{theo}

\zerarcounters
\section{Numerical results}
\label{sec:4}

Numerically we may illustrate the scenarios considered in the previous sections
as well as provide some insight into cases which are not covered by the analysis.
We study the system \eqref{eq:1.1} with $f(t) = \cos t $. In particular we consider three situations,
$|\mu| < 1$, $|\mu| > 1$ and $|\mu| = 1$, with both high and low frequency forcing.

To check the stability of the dynamics we take 10 000 pseudo-random initial conditions
within a square $[-2,2]\times[-2,2]$ from the phase plane $(x,y)$.
The chosen numerical integration method is a St\"{o}rmer-Verlet scheme with variable step size.
The St\"{o}rmer-Verlet method is a second order symplectic scheme, details of which may be found in \cite{LR}.

After an initial transient period, the trajectories are checked to ascertain how their asymptotic behaviour
has changed with respect to the trajectories of the unperturbed system with the same initial conditions.
If most of the orbits have remained close to the corresponding orbits of the unperturbed system,
then we say that the system is ``stable".
This is expected to occur when the system is well within the KAM regime: the majority of the unperturbed tori
persist slightly deformed, so that every orbit either lies on a torus or is trapped between two surviving tori.
However it is possible that the trajectories of the perturbed system do not remain close, but are still bounded.
This can happen as we are moving out of the KAM regime: most of the tori are destroyed,
with a few of them still existing and undergoing much larger deformations.
We refer to such a case by saying that the system is ``bounded". Numerically it is difficult to classify trajectories
as unbounded, as a trajectory which appears unbounded may be bounded within a very large region.
Therefore, pragmatically, we class the trajectories as unbounded once their amplitude exceeds 30
in either the $x$ or $y$ direction, and class the system as ``unbounded" if any trajectory is found to be so.
When this happens, nearly all (if not all) KAM tori are expected to be destroyed,
at least in the region investigated, otherwise any of them would confine the orbits inside.
We note, however, that, even though KAM theory no longer applies in this case,
one can still have invariant curves of a different kind, such as cantori (see for instance \cite{R}),
so that it is quite possible for some trajectories to be unbounded whilst others remain bounded
(and possibly not too far from those of the unperturbed system).

First we consider the scenario where $|\mu| < 1$, in particular we choose $\mu = 0.8$, 0.9 and
0.95; some numerical results are shown in Table \ref{TablemuLessThan1}.
For $\omega$ sufficiently large or small,
as proved in the previous sections, the dynamics are stable,
whilst between these two extremes the system loses stability.

\vspace{.4cm}
\begin{table}[H]
\centering
\small{
\begin{tabular}{|cccc|}
\hline\hline
\multicolumn{1}{|c}{$\omega$} & \multicolumn{1}{c}{$\mu = 0.8$} & \multicolumn{1}{c}{$\mu = 0.9$}& \multicolumn{1}{c|}{$\mu = 0.95$} \\
\hline\hline
0.05 & Stable & Stable & Stable \\
0.1 & Stable & Stable & Bounded \\
0.2 & Stable & Bounded & Bounded \\
0.3 & Bounded & Bounded & Unbounded \\
0.4 & Bounded & Unbounded & Unbounded \\
0.5 & Bounded & Unbounded & Unbounded \\
0.7 & Bounded & Unbounded & Unbounded \\
0.8 & Unbounded & Unbounded & Unbounded \\ \hline
12.0 & Unbounded & Unbounded & Unbounded \\
13.0 & Stable & Unbounded & Unbounded \\
14.0 & Stable & Stable & Stable \\ \hline
\end{tabular}
}
\caption{\small{Stability and boundedness results for choices of $|\mu| < 1$ and various values of $\omega$.
We class the system as ``stable'' if most of the orbits remain close to those of the unperturbed system, ``bounded'' if
none of the orbits exceed 30 in either direction and ``unbounded'' if it is not bounded,
i.e. at least one orbit is unbounded; see text for details.}}
\label{TablemuLessThan1}
\end{table}

\begin{figure}[H]
\centering
\subfloat[]{\includegraphics[width=0.36\textwidth]{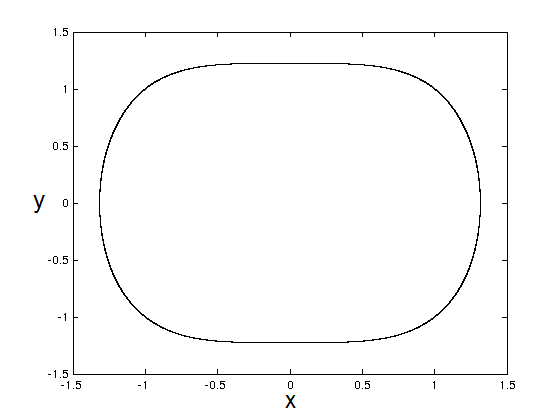}}
\subfloat[]{\includegraphics[width=0.36\textwidth]{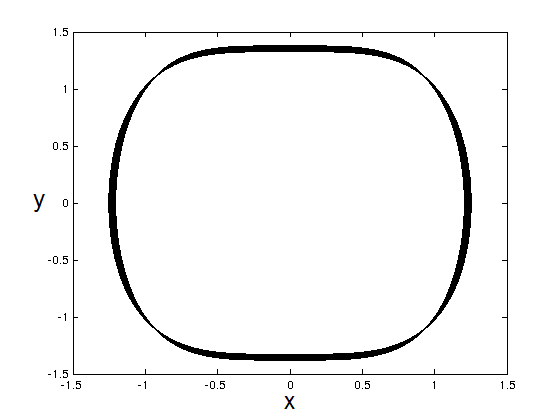}}\\
\subfloat[]{\includegraphics[width=0.36\textwidth]{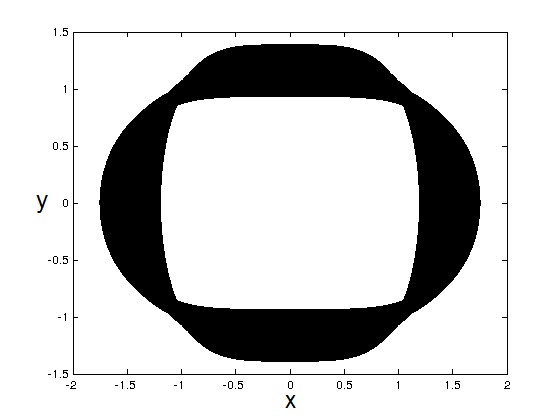}}
\subfloat[]{\includegraphics[width=0.36\textwidth]{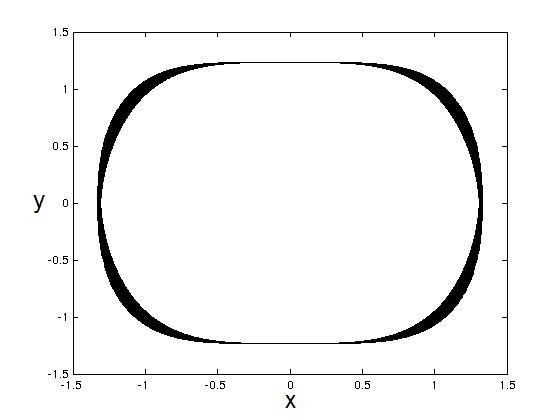}}
\caption{\small{Orbits for the system \eqref{eq:1.2} with initial conditions
$(x,y) = (1,1)$. In Figure (a) $\mu = 0$. In Figures (b), (c) and (d) $\mu = 0.8$
and $\omega = 0.0001$, 0.2 and 14, respectively.
}}
\label{FiguremuLessThan1}
\end{figure}

In Figure \ref{FiguremuLessThan1}
we show some example orbits corresponding to the initial conditions $(x,y) = (1,1)$
with $\mu = 0$ in Figure \ref{FiguremuLessThan1}(a) and $\mu = 0.8$ in (b), (c) and (d).
It may be seen that the perturbed orbits remain close to the unperturbed system
for suitable $\omega$. With $\mu$ getting close to $1$, one has to take $\om$ increasingly small
or increasingly large for the system to be stable.

\vspace{-.1cm}

\vspace{-.5cm}
\begin{figure}[H]
\centering
\subfloat[]{\includegraphics[width=0.4\textwidth]{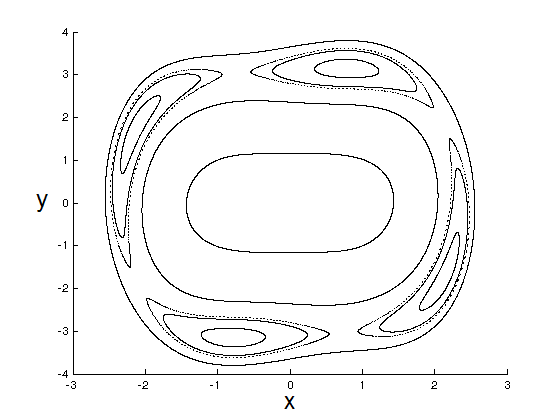}}
\subfloat[]{\includegraphics[width=0.4\textwidth]{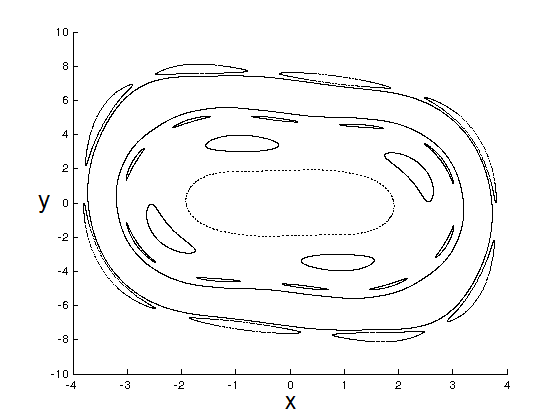}}
\caption{\small{Poincar\'{e} maps for the system \eqref{eq:1.1} with $\mu = 0.8$
and various initial conditions, showing the break up of the KAM tori.
The frequency is $\omega = 0.2$ and 0.4 in (a) and (b), respectively.
}}
\label{PoincareFigure}
\end{figure}

The results in Table \ref{TablemuLessThan1} show that for $\mu = 0.8$ the system becomes unstable
when $\omega > 0.2$. We see in Figure \ref{PoincareFigure}(b) and (c) that as $\omega$ increases from
0.2 to 0.4, many of the tori are broken and cantori appear, separating the few persisting KAM tori and
creating larger and larger gaps. However, as quite a few tori still persist in the region investigated,
the orbits remain bounded. Increasing $\omega$ further causes more of the KAM tori to break up
and for $\omega = 0.8$ the orbits are no longer bounded in the region $[-30,30]\times[-30,30]$.
Analogous considerations apply to the other values of $\mu$: the closer the amplitude $\mu$ is to 1,
the more extreme $\omega$ must be for the system to be stable.

\vspace{.4cm}
\begin{table}[H]
\centering
\small{
\begin{tabular}{|ccccc|}
\hline\hline
\multicolumn{1}{|c}{$\omega$} & \multicolumn{1}{c}{$\mu = 1$} & \multicolumn{1}{c}{$\mu = 1.2$} & \multicolumn{1}{c}{$\mu = 2$} & \multicolumn{1}{c|}{$\mu = 5$}\\
\hline\hline
0.0001 & Stable & Stable & Stable & Stable \\
0.0002 & Stable & Stable & Bounded & Unbounded \\
0.0003 & Bounded & Unbounded & Unbounded & Unbounded \\
0.0004 & Bounded & Unbounded & Unbounded & Unbounded \\
0.0005 & Bounded & Unbounded & Unbounded & Unbounded \\
0.0010 & Unbounded & Unbounded & Unbounded & Unbounded \\ \hline
13 & Unbounded & Unbounded & Unbounded & Unbounded \\
14 & Stable & Stable & Unbounded & Unbounded \\
16 & Stable & Stable & Unbounded & Unbounded \\
18 & Stable & Stable & Stable & Unbounded \\
20 & Stable & Stable & Stable & Stable \\ \hline
\end{tabular}
}\caption{\small{Stability and boundedness results for choices of $|\mu| \ge 1$ and various values of
$\omega$. The system is classed as either ``stable" or ``bounded" or ``unbounded"
as explained in the caption of Figure \ref{TablemuLessThan1}.}}
\label{TablemuMoreThan1}
\end{table}

\begin{figure}[H]
\centering
\subfloat[]{\includegraphics[width=0.33\textwidth]{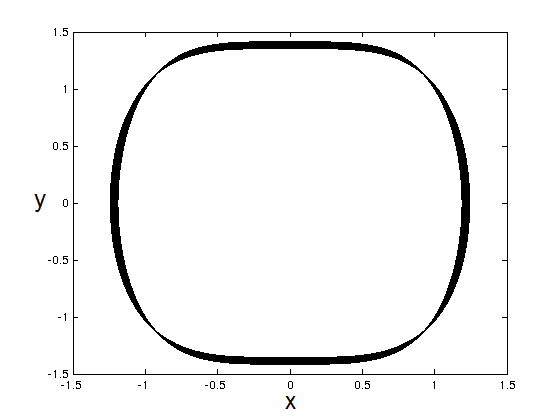}}
\subfloat[]{\includegraphics[width=0.33\textwidth]{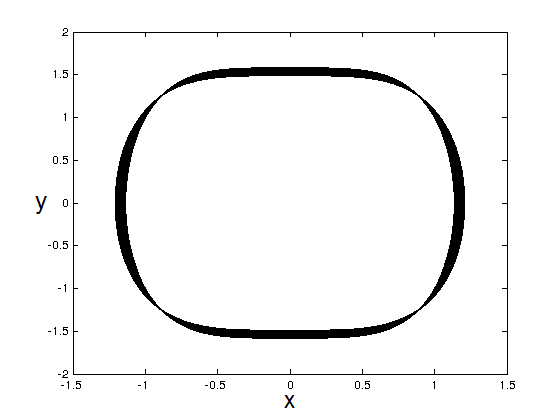}}
\subfloat[]{\includegraphics[width=0.33\textwidth]{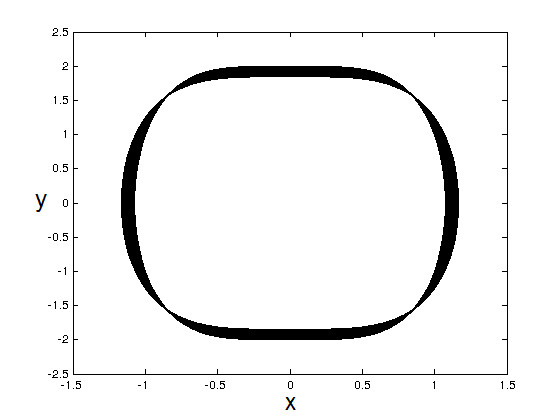}}\\
\subfloat[]{\includegraphics[width=0.33\textwidth]{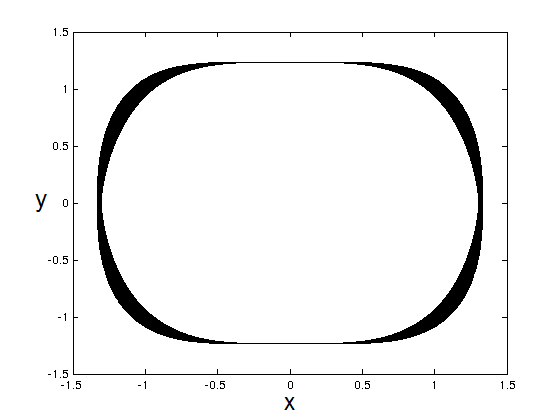}}
\subfloat[]{\includegraphics[width=0.33\textwidth]{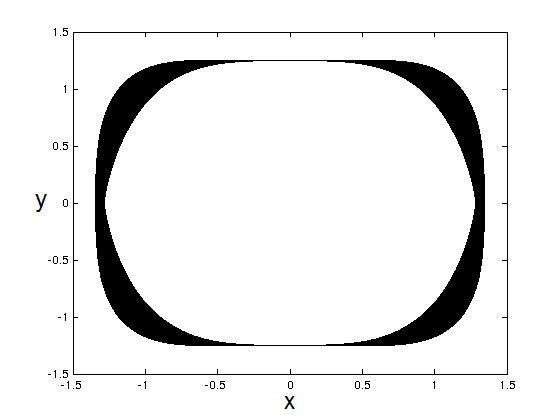}}
\subfloat[]{\includegraphics[width=0.33\textwidth]{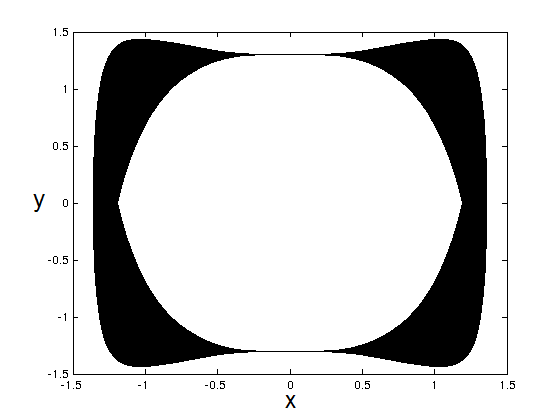}}
\caption{\small{Orbits for the system \eqref{eq:1.2} with initial conditions
$(x,y) = (1,1)$. In Figures (a),(b) and (c) $\omega = 0.0001$ and in Figures
(d), (e) and (f) $\omega = 14$. The amplitude is $\mu = 1$, 2 and 5 in Figures (a, d), (b, e)
and (c, f), respectively.
}}
\label{FiguremuMoreThan1}
\end{figure}

For $|\mu| \ge 1$ the analysis in the previous sections can only be applied when the system undergoes high frequency forcing.
Numerically we find that the system is also stable with low frequency forcing, however $\omega$ must be taken considerably
smaller than the cases where $|\mu| < 1$. This is not true for high frequency forcing, where similar orders of $\omega$
(compared with the cases where $|\mu|<1$) are sufficient for the dynamics to become stable. Some numerical
results are presented in Table \ref{TablemuMoreThan1}.
In Figure \ref{FiguremuMoreThan1} we show some example orbits with $\mu \ge 1$.

Similarly to the case where $|\mu| < 1$ it is evident
that the perturbed orbits remain close to that of the unperturbed system provided $\omega$ is
suitably chosen. In Figure \ref{FiguremuMoreThan1}(f) we see that, although the system
is classed as unbounded for $\omega = 14$ and $\mu = 5$, it is still possible to find initial
conditions for which the orbit remains close to the unperturbed system.

\zerarcounters
\section{Conclusions}
\label{sec:5}

More generally the Hamiltonian of a forced cubic oscillator is
\begin{equation} \label{eq:4.1}
H(y,x,t)= \frac{1}{2} m y^{2} + \left( a + \m f(\omega t) \right) \frac{x^{4} }{4} .
\end{equation}
However the writing \eqref{eq:1.1} is not restrictive, since we can reduce \eqref{eq:4.1}
to that form by rescaling both variables $x$ and $y$ and redefining the parameter $\mu$.

\vspace{-.4cm}

Other generalisations can be easily envisaged. For instance any potential
$V(x)$ yielding closed orbits in a region encircling the origin can be considered.
In particular one can take a potential $x^{2n}/2n$ instead of $x^4/4$: the unperturbed
system is still integrable, so that the analysis of Section \ref{sec:2}
applies immediately. Also the discussion in Section \ref{sec:3} can be easily
adapted to cover such a case; we refer to \cite{DZ2,LL} for details.
Also, less regularity can be required for the driving force.
Finally, one could consider quasi-periodically forced systems, as in \cite{LZ},
in the case in which all components of the frequency vectors are small or large.

Coming back to our system \eqref{eq:1.1}, the condition that $f$ has zero average
could be relaxed. In fact, for Theorem \ref{thm:2} to hold, what we really need is that
$1 + \mu \langle f(\cdot) \rangle >0$, so that the averaged system is integrable.
On the other hand, Theorem \ref{thm:3} requires $1+ \mu f(\cdot)>0$.
Thus, if $\langle f(\cdot)\rangle=0$ and $\|f\|_{\io}=1$, this excludes the case $\mu=1$.
However, if the average of $1+\mu f(\omega t)$ is positive, one can argue that the potential remains
positive for most of time, so one can conjecture that the condition $1+\mu \langle f(\cdot) \rangle>0$
might be sufficient in the low frequency regime as well. As discussed in Section \ref{sec:4},
we have found numerically that also in such a case the orbits are stable if the forcing frequency
is sufficiently low. It would be interesting to investigate the issue in more detail by means of analysis.

\end{document}